\documentclass[12pt,twoside]{article}

\usepackage{amsmath}
\usepackage{amsthm}
\usepackage{latexsym}
\usepackage{amssymb}

\newcommand{\qkappa}{{\mathbb Q}(\kappa)}
\newcommand{\QQ}{{\mathbb Q}}
\newcommand{\rest}{\!\upharpoonright\!}

\newcommand{\bAA}{{\cal A}}
\newcommand{\HH}{{\cal H}}

\newcommand{\deq}{\buildrel{\rm def}\over =}
\newcommand{\hal}{{\empty^{\kappa}}{\cal H}^{\rho}_{\alpha}}
\newcommand{\hnal}{{\empty^{\kappa}}{\cal H}^{\rho}}
\newcommand{\eop}{\hfill$\Box$ \newline}

\theoremstyle{definition}
\newtheorem{Theorem}{Theorem}[section]
\newtheorem{Main Theorem}[Theorem]{Main Theorem}
\newtheorem{Lemma}[Theorem]{Lemma}
\newtheorem{Subclaim}[Theorem]{Subclaim}
\newtheorem{Claim}[Theorem]{Claim}

\newtheorem{Definition}[Theorem]{Definition}
\newtheorem{Notation}[Theorem]{Notation}
\newtheorem{Remark}[Theorem]{Remark}

\newtheorem{Question}[Theorem]{Question}
\newtheorem{Corollary}[Theorem]{Corollary}
\newtheorem{Observation}[Theorem]{Observation}

\newenvironment{Proof of the Subfact}{\noindent{\bf Proof of the Subfact.}}{\eop}
\newenvironment{Proof of the Theorem}{\noindent{\bf Proof of the Theorem.}}{\eop}
\newenvironment{Proof of the Conclusion}{\noindent{\bf Proof of the Conclusion.}}{\eop}
\newenvironment{Proof of the Observation}{\noindent{\bf Proof of the Observation.}}{\eop}
\newenvironment{Proof of the Fact}{\noindent{\bf Proof of the Fact.}}{\eop}
\newenvironment{Proof of the Lemma}{\noindent{\bf Proof of the Lemma.}}{\eop}
\newenvironment{Proof of the Claim}{\noindent{\bf Proof of the Claim.}}{\eop}
\newenvironment{Proof of the Subclaim}{\noindent{\bf Proof of the Subclaim.}}{\eop}
\newenvironment{Proof of the Main Claim}{\noindent{\bf Proof of the Main Claim.}}{\eop}
\newenvironment{Proof}{\noindent{\bf Proof.}}{\eop}

\title{An Expansion of a Poset Hierarchy} 
\author{Mirna D{\v z}amonja and Katherine Thompson}

\begin{document}
\maketitle
\begin{abstract} 
This article extends a paper of Abraham and Bonnet which 
generalised the famous Hausdorff characterisation of the class of scattered
linear orders. Abraham and Bonnet gave a poset hierarchy that characterised the class of
scattered posets which do not have infinite antichains (abbreviated
FAC for finite antichain condition). An antichain here is taken in the sense
of incomparability.
We define a larger poset hierarchy than that of Abraham and Bonnet,
to include a broader class of ``scattered'' posets that we
call $\kappa$-scattered. These posets cannot embed any order such that
for every two subsets of size $ < \kappa$, one being 
strictly less than the other, there is an element in
between. If a linear
order has this property and has size $\kappa$ we call
this set $\qkappa$. Such a set only exists when
$\kappa^{<\kappa}=\kappa$. Partial orders with
the property that for 
every $a<b$ the set $\{x:\,a<x<b\}$ has size $\geq \kappa$ are called weakly $\kappa$-dense,
and partial orders that do not have a weakly $\kappa$-dense subset are called
strongly $\kappa$-scattered.
We prove that our hierarchy includes all strongly $\kappa$-scattered
FAC posets, and that the hierarchy
is included in the class of all FAC $\kappa$-scattered posets.
In addition, we prove that our hierarchy is in fact the closure of the class
of all $\kappa$-well-founded linear orders under inversions, lexicographic sums and
FAC weakenings.
For $\kappa=\aleph_0$ our hierarchy agrees with the one from the Abraham-Bonnet
theorem.
\footnote{The authors warmly thank Uri Abraham for his many useful suggestions
and comments.
This research was started in the second author's undergraduate
thesis supervised by
the first author
at the University of Wisconsin and was partly funded by the
University League
Scholarship provided
by the College of Letters and Science in the University of Wisconsin.
Mirna D\v zamonja thanks EPSRC for their support on an EPSRC Advanced Fellowship.}
\end{abstract}
{\em Keywords: Set theory, ordered sets, $\kappa$-dense, $\eta_{\alpha}$-orderings\\
MSC 2000 Classification numbers: 03E04, 06A05, 06A06}

\section{Introduction} Recall that a scattered order is one which does
not embed the rationals.  
Hausdorff (\cite{Ha}, or see \cite{haus})
proved that the class of scattered linear orders is the
least family of linear orderings which includes the ordinals and is closed under
lexicographic sums
and inversions. The paper \cite{AB}
by Abraham and Bonnet proved that the class of scattered posets
satisfying FAC (the finite antichain condition) 
is the least family of posets satisfying FAC which includes the
well-founded posets and is closed
under inversions, lexicographic sums and augmentations. 

There are several routes for expansion on these results which centre
around a generalisation of the concept of scattered to higher cardinalities. To this effect,
one would consider a $\kappa$-scattered poset (or linear order) to be
one which does not embed a $\kappa$-dense set. There are two
definitions that one could give to a $\kappa$-dense set. The first was
introduced by Hausdorff in 1908 as an $\eta_{\alpha}$-ordering for
$\kappa = \aleph_\alpha$.  This
is an order such that between any two subsets of size $ < \kappa$, one being
strictly less than the other, there is an element in between.  Orders
with this property are here called strongly $\kappa$-dense. When 
an $\eta_{\alpha}$-ordering is linear and also has size $\kappa$, we call it $\qkappa$. The
other definition of $\kappa$-dense is a strictly weaker one in which between
every two elements there is a subset of size $\kappa$. We call this notion
weakly $\kappa$-dense.
Using either of the two definitions for $\kappa$-scattered, which
we will call weakly $\kappa$-scattered (not embedding a strongly $\kappa$-dense set)
and strongly $\kappa$-scattered (not embedding even a weakly $\kappa$-dense set),
we can
attempt to expand the characterisation results on linear orders or FAC
posets. Note that the class of strongly $\kappa$-scattered is included in the
class of weakly $\kappa$-scattered orders.

This paper builds on \cite{AB} and extends its results. As in
\cite{AB}, a class of posets is built in a hierarchical way so that 
for any regular $\kappa$
we have that ${}^\kappa\HH$ is the least family
of posets satisfying FAC which includes the  $\kappa$-well founded
posets and is closed under inversions, lexicographic 
sums and augmentations. We then close this class under
FAC weakenings (the dual notion to augmentations, but retaining the
FAC) to obtain the class ${}^\kappa\HH^\ast$.  
We prove that class ${}^\kappa\HH^\ast$ contains
all strongly $\kappa$-scattered posets and is contained in the class of all
weakly $\kappa$-scattered FAC posets. 
For $\kappa=\aleph_0$ where the two notions
of scattered agree the two hierarchies agree, as follows
by the Abraham-Bonnet theorem, and are both equal to the class of
FAC scattered posets.

It is also shown that the class ${}^\kappa\HH^\ast$ can be constructed
in a simpler way, by starting with the $\kappa$-well founded linear
orders and closing under inversions, lexicographic sums and FAC
weakenings. It is proved that this is exactly the same class as the
one constructed by posets.

A reader familiar with \cite{AB}
may at this point wonder why it is that for $\kappa>\aleph_0$ we do not obtain the complete
analogue of the Abraham-Bonnet theorem. There are two main difficulties, apart from the
fact that the notions of weakly and strongly $\kappa$-scattered for $\kappa>\aleph_0$
are distinct, as opposed to what happens at $\kappa=\aleph_0$.
The first one is that it is not necessarily the case that if all augmentations of a poset 
are weakly or strongly $\kappa$-scattered then
the poset has the FAC, or at least we have not been able to prove this.
The other difficulty is
that we do not know how to prove that FAC posets which
are not in the hierarchy defined
above actually embed a strongly $\kappa$-dense set, although we can prove that 
they have an augmentation that embeds a weakly $\kappa$-dense subset.

It remains unknown whether every weakly $\kappa$-scattered poset is in the
hierarchy ${}^\kappa\HH^\ast$ or if
${}^\kappa\HH^\ast$ and ${}^\kappa\HH$ are in general equal. However, ${}^\kappa\HH$ does
contain examples of weakly $\kappa$-dense posets (as we will show in
the final section), so it cannot be the case that ${}^\kappa\HH$ only contains strongly
$\kappa$-scattered posets. 

\section{On $\kappa$-scattered Posets}

We start by explaining how Abraham and Bonnet's theorem extends Hausdorff's theorem.
We first need several definitions. In this paper, we use `order' to denote a `partial
order', and whenever we deal with linear orders we specify this.

\begin{Definition} 
\begin{description}
\item{(1)} A (partial) order $P$ {\em embeds} an order $Q$ iff there is an order
preserving one-to-one function from $Q$ to $P$.

\item{(2)} An order is said to be {\em scattered} iff it does not embed the rationals,
${\mathbb Q}$, with their usual ordering.

\item{(3)} If $(I, \le_I)$ is a partial order and $\bar{P}=\langle (P_i,
  \le_i):\,i\in I\rangle$ is a sequence of  partial orders, 
the {\em lexicographic sum} of $\bar{P}$ is the order whose
universe is $\bigcup_{i\in I}P_i$, ordered by letting
$p\le q$ if and only if there exists $i\in I$ such that $p\le_i q$ or there exists 
$i<_I j$ with $p\in P_i$ and  $q\in P_j$.
\item{(4)} A poset $(P',\le')$ is an {\em augmentation} of a poset $(P, \le)$ if and only if
$P = P'$, and for all $p,q\in P$, if $p\le q$ then $p\le' q$. We also say that
$P$ is a {\em weakening} of $P'$.
\item{(5)} If $P$ is a subposet of $Q$ in which all
relevant $Q$-relations are kept then $P$ is said to be a {\em restriction} of 
$Q$ to $P$ (written $Q \rest P$). In particular, for any $p_1, p_2 \in P$ we have $p_1 <_P p_2$ 
iff $p_1 <_Q p_2$. 
\item{(6)} We say that a poset $P$ is {\em $\kappa$-well founded} if and only if
$P$ does not have any decreasing sequences of size $\kappa$. 
\end{description}
\end{Definition}

Note that the notion of $\kappa$-well founded was introduced by Zaguia
as the extraction property for $\kappa$, see \cite[\S 4.11.3]{fraisse}. 

The relevant theorem  of Hausdorff in \cite{Ha} (see also \cite{haus})
states that the scattered linear orders are exactly the
closure of the class of all well orderings under inversions and
lexicographic sums.

\begin{Notation} \begin{enumerate} \item Let  $p \bot q$ denote that $p$ is incomparable
with $q$.
\item We will say {\em antichain} when we mean an incomparable antichain, that is, a 
subset whose elements are pairwise incomparable.
\item Let FAC denote the finite antichain condition, so in an FAC poset all
antichains are finite.
\item Similarly, let $\kappa$-AC denote the $\kappa$-antichain condition. That is, a poset 
has the $\kappa$-AC if and only if it does not have an antichain of size $\kappa$. 
\item  If $(P, \leq)$ is a poset and  $S, T \subseteq P$, we write
$S\le T$ iff for all $s\in S$ and $t\in T$, we have $s\le t$.
\end{enumerate} \end{Notation}

Clearly, linear orders are a special case of FAC posets. Abraham and
Bonnet proved that the class of scattered FAC posets is the closure of
the class of well founded FAC posets under augmentations, inverses and lexicographic sums. 
Some of the main tools they used were Hessenberg based operations on ordinals and the
notion of the antichain rank $\rho$. We shall not need to reintroduce
Hessenberg operation, and as for
the antichain rank, basically it is a function that determines the length of the set of antichains 
in any given FAC poset. We include the definition given in \cite{AB} here; it will
be needed in \S\ref{mainsection}.

\begin{Definition}\label{ranks} For any FAC poset $P$, let $(\bAA(P), \supset)$ be the poset of 
all non-empty antichains of $P$ under inverse inclusion. Since this is a well-founded
poset, we can define the usual rank function on it which we will call the antichain
rank of $P$ and denote by rk$_{\bAA}(P) \deq \mbox{ rk}(\bAA(P))$.  
\end{Definition}

Hausdorff's theorem is in fact the restriction of the Abraham-Bonnet theorem
to antichain rank 1. 
Let us now go on to
define what we mean by weakly $\kappa$-scattered, by first defining the dual notion, 
strongly $\kappa$-dense. 

\begin{Definition}\label{exist-star} (1) For
a cardinal $\kappa \geq \aleph_0$ we say an order $(P, <^*)$ is strongly
$\kappa$-dense iff

\hspace{2cm} $(\forall S,T\in [P]^{<\kappa})\,[S<^\ast T\implies (\exists x) S<^\ast x <^\ast T],$
\hfill $(*)^{\kappa}$. \\

(2) We denote by $\qkappa$ a strongly $\kappa$-dense linear order of
    size $\kappa$ whenever this set exists and is unique up to isomorphism.

\end{Definition}

A linear order which satisfies $(*)^{\kappa}$ is also known as an
$\eta_{\alpha}$-ordering for $\kappa = \aleph_{\alpha}$. Hausdorff
proved in \cite{Ha} that such an ordering exists for all regular
cardinals $\kappa$. However, it can only be shown that these sets can have
size $\kappa$ when $\kappa$ satisfies this stronger property, $\kappa
= \kappa^{< \kappa}$.
 
 We know that the countable version of this set exists,
namely the rationals satisfy this for $\kappa=\aleph_0$. 
It follows from Shelah's work on the existence
of saturated models for unstable theories (see \cite{Shbookc}, Theorem
VIII 4.7)
that $\qkappa$ exists iff $\kappa^{<\kappa}=\kappa$. The specific
instance of this result for a dense linear order with no endpoints
is well-known.

Sierpinski showed that for general $\kappa$ satisfying $\kappa^{< \kappa} = \kappa$
the order $\qkappa$ may be constructed by induction (see
\cite{linorders} for details). The same proof also gives a more
general construction of an order  of size $\lambda$ which
satisfies $(*)^{\kappa}$ where $\kappa$ is a regular
cardinal and $\kappa^{<\kappa} = \lambda$. 


                                                                                 
                                                                                 

The obvious way to generalise the notion of scattered would be to say
that an order is $\kappa$-scattered iff it does not embed the unique
linear order $\qkappa$. However, since this set only exists given some
strong cardinal arithmetic assumptions, Stevo Todor{\v c}evi{\' c} suggested
that it is more natural to 
say that an order is $\kappa$-scattered iff it does not embed a
strongly $\kappa$-dense set of any size. In this way, the notion makes
sense whenever $\kappa$ is a regular cardinal. 

The next claim shows that whenever $\qkappa$ exists, the properties of
being strongly $\kappa$-dense and embedding $\qkappa$ are equivalent.

\begin{Claim} Suppose $\qkappa$ exists and $P$ is a strongly
  $\kappa$-dense order. Then there is $Q \subseteq P$ such that $Q$ is
  isomorphic to $\qkappa$.
\end{Claim}

\begin{Proof of the Claim} As shown above, $\qkappa$ exists whenever
  $\kappa = \kappa^{< \kappa}$. Since $P$ is strongly $\kappa$-dense,
  it follows that $|P| \geq \kappa$. Let $Q_0$ be any subset of $P$ of
  size $\kappa$. In particular, $|Q_0|^{< \kappa} = \kappa$.

Let $\{A_{\alpha}, B_{\alpha} : \alpha < \kappa\}$ enumerate all $A, B
\in [Q_0]^{< \kappa}$ such that $A < B$. For each $\alpha$, there
exists $x_{\alpha} \in P$ such that $A_{\alpha} < x_{\alpha} <
B_{\alpha}$ by the fact that $P$ is strongly $\kappa$-dense. 
Let $Q_1 = Q_0 \cup \{x_{\alpha} : \alpha < \kappa\}$.

Repeat this process inductively, creating $Q_\beta$ at each step
$\beta < \kappa$ taking the union at limit stages. Let $Q =
\bigcup_{\beta < \kappa} Q_\beta$.

We will show that the $Q$ we have constructed is isomorphic to
$\qkappa$. First at each successor stage $\beta + 1 < \kappa$ we have added only
$\kappa$ many elements to $Q_\beta$ so $|Q_{\beta+ 1}| = \kappa$. Note
that no new elements are added at limit stages. Therefore, $Q$ is
the union of $\kappa$ many sets of size $\kappa$ so itself has size
$\kappa$. 

For $A, B \in Q$ such that $|A|, |B| < \kappa$ and $A < B$
there exists $\beta < \kappa$ such that $A, B \subseteq Q_{\beta}$ as
$\kappa$ is regular. By the construction, there exists $x \in Q_{\beta
  + 1}$ such that $A < x < B$. Thus, $Q$ is strongly $\kappa$-dense. 
\end{Proof of the Claim}

Thus for example, the notion of strongly $\kappa$-dense
for $\kappa=\aleph_0$ agrees with the definition
of dense as $\QQ(\aleph_0) = \QQ$ and so an order which embeds a
strongly $\aleph_0$-dense set also embeds the rationals.

The following fact about strongly $\kappa$-dense sets will be useful
to us. 

\begin{Lemma}\label{q-dec} Any strongly $\kappa$-dense set necessarily
  has a $\kappa$-decreasing sequence 
and a $\kappa$-increasing sequence.
\end{Lemma}
\begin{Proof of the Lemma} Let $P$ be a strongly $\kappa$-dense set.
  We prove this for a $\kappa$-decreasing sequence, 
the proof for a $\kappa$-increasing sequence is similar.
By induction on $\alpha < \kappa$, we construct 
$A = \{q_{\alpha} : \alpha < \kappa\}$ with $q_{\alpha} >^* q_{\beta}$ if $\alpha < \beta$. Let
$q_0$ be any element of $P$. 

For any $0 < \alpha < \kappa$, assume $q_{\beta}$ is defined for all
$\beta < \alpha$. Let $T = \{q_{\beta} : \beta < \alpha\}$ and $S = \emptyset$.
Choose $q_{\alpha}$ to be any $q$ such that  $S <^\ast q <^\ast T$. 
\end{Proof of the Lemma}

Now that we have proved the relevant properties of strongly
$\kappa$-dense sets, we may turn to their opposite, the idea of weakly
$\kappa$-scattered sets. 

\begin{Definition}\label{scatt} Suppose that $\kappa$ is a regular cardinal. We say
that a partial order is weakly $\kappa$-scattered if and only if it does not embed
any strongly $\kappa$-dense set. We may omit the adjective `weakly' when discussing this notion.
\end{Definition}

Hence for $\kappa$ as in Definition \ref{scatt},
all orders of size
$<\kappa$, in particular finite orders, are $\kappa$-scattered.
If $\kappa>\aleph_0$, then there are orders which are $\kappa$-scattered
and not scattered, for example the rationals. Similarly, if
$\kappa_1>\kappa_2$ both satisfy the assumptions of Definition \ref{scatt},
then there are orders which are $\kappa_1$-scattered without being
$\kappa_2$-scattered. In the other direction, every $\kappa_2$-scattered
order is $\kappa_1$-scattered, as 
we can see that in these circumstances ${\mathbb Q}(\kappa_2)$
embeds into ${\mathbb Q}(\kappa_1)$ whenever these sets exist.

Our aim is to consider the Abraham-Bonnet theorem for $\kappa$-scattered
partial orders which satisfy FAC for regular cardinals
$\kappa$ with $\kappa \geq \aleph_0$.
We shall start by showing that
strongly $\kappa$-dense sets have a property which might seem stronger than
$(\ast)^\kappa$, but as the claim shows, is actually equivalent to  
it. The proof is similar to that of Lemma \ref{q-dec}.

\begin{Claim}\label{star} Suppose $P$ is a poset satisfying $(\ast)^\kappa$.
Then for all $S, T \subseteq P$ with $|S|, |T| < \kappa$ and
$S<_PT$ we have $|\{q \in P : S <_P q <_P T\}| \geq \kappa$. Moreover, if
$\kappa = \kappa^{< \kappa}$ and
$P$ is a linear order then $P$ restricted to the suborder
$Q = \{q \in P :$ $S <_P q <_P T\}$ is isomorphic to $\qkappa$.
\end{Claim}

\begin{Proof of the Claim} Let $S, T$ satisfy the assumptions of the claim.  
By induction on $\alpha < \kappa$, we construct $\{q_\alpha : \alpha < \kappa\}$ such that
$\beta < \alpha$ implies $q_{\beta} <_P q_\alpha$ and $S <_P q_\alpha <_P T$. Let $q_0$ be any
such $q$,
which exists by the assumption $(\ast)^\kappa$.  
Assume we are given $\{q_{\beta} : \beta < \alpha\}$, by the induction
hypothesis.  Apply $(\ast)^\kappa$ to
$T^* = T \cup \{q_{\beta} : \beta < \alpha\}$ and $S$, noticing that $|T^*| < \kappa$ and
$S<_P T^\ast$ to obtain $q_{\alpha}$.  

If $P$ satisfies $(\ast)^\kappa$, then the suborder $Q$ as defined in the claim
must also satisfy  $(\ast)^\kappa$. Therefore if $\kappa = \kappa^{<
  \kappa}$ then $Q$ is a linear order of size $\kappa$ which satisfies
$(\ast)^\kappa$ and hence isomorphic to $\qkappa$.  
\end{Proof of the Claim}

In our main
result we shall use a weaker notion of $\kappa$-density as well, so we define it here. 

\begin{Definition} If a linear order $L$ satisfies the property
\[ |L|\ge 2
\mbox{ and }(\forall s,t\in L)[s< t\implies |\{x:\,s< x<  t\}|\ge\kappa]
\]
then we say that $L$ is weakly $\kappa$-dense. We may omit the adjective `weakly'
when discussing this notion.
(The first clause of the property is included to avoid trivialities.)

An order that does not embed a weakly $\kappa$-dense order is called 
strongly $\kappa$-scattered.
\end{Definition}

For $\kappa>\aleph_0$
it is easy to construct an example of a $\kappa$-dense linear order that is not
strongly $\kappa$-dense, moreover there are $\kappa$-dense linear orders
that are $\kappa$-scattered and ones that do not even have a decreasing
$\kappa$-sequence. See \S\ref{examples}. Note that at $\kappa =
\aleph_0$, the two definitions of $\kappa$-density agree. 

If an order is not $\kappa$-scattered for $\kappa = \kappa^{<
  \kappa}$  then it embeds a copy of $\qkappa$ so clearly
it has a suborder of size $\kappa$ that is not $\kappa$-scattered. For future
purposes we note that a similar statement is true about orders that are
strongly $\kappa$-scattered for any cardinal $\kappa$.

\begin{Claim}\label{downtokappa} Suppose that $P$ is an order that is not
strongly $\kappa$-scattered for $\kappa  = \kappa^{<
  \kappa}$. Then $P$ has a suborder of size $\kappa$ that is
not strongly $\kappa$-scattered.
\end{Claim}

\begin{Proof of the Claim} We shall define a suborder $Q=\bigcup_{n<\omega} Q_n$ of $P$
by defining $Q_n$ by induction on $n$. Let $Q_0$ be any two-element linear
suborder of $P$, which exists by the definition of weak $\kappa$-density.
Given $Q_n$ of size $\le\kappa$ let us choose for any $a<_{Q_n}b$ a set 
$S^n_{a,b}\subseteq P$ of size $\kappa$ such that $a< S^n_{a,b}< b$ and let
$Q_{n+1}=\bigcup_{a<_{Q_n}b}S^n_{a,b}\cup Q_n$. It is easy to see that $Q$
is as required.
\end{Proof of the Claim}

It is a well-known theorem of Bonnet and Pouzet
and independently
Galvin and McKenzie (see \cite{BP}), that every scattered partial order has
a scattered linear extension. 
An important ingredient in the Abraham-Bonnet theorem is a lemma which
says that an FAC partial order is scattered iff all its augmentations
are scattered. 
In our situation, we shall not be able to get such a neat
equivalence, but a chain of implications instead. To prove the mentioned
equivalence, Abraham and Bonnet use a particular claim which
relies heavily on the fact that the
lexicographic sum along ${0<1}$ (i.e. the union) of two scattered partial
orders is still scattered. In our circumstances we need
the following claim. 

\begin{Claim}\label{mirna} Assume that $\kappa$ is a regular cardinal.
Suppose that $(D, <_D)$ is a poset of size $\kappa$ and 
$D=\bigcup_{i<i^\ast} D_i$ for some $i^\ast <\kappa$, and each
$D_i$ is $\kappa$-scattered. Then $D$ is $\kappa$-scattered.
\end{Claim}

\begin{Proof of the Claim} The proof is by induction on $i^\ast < \kappa$. For 
$i<i^\ast$ let $R_i\deq\bigcup_{j<i} D_j$.
By the induction hypothesis, each $R_i$ is $\kappa$-scattered, and by
definition, we have $D=\bigcup_{i<i^\ast}R_i$, while
$j<i<i^\ast$ implies $R_j\subseteq R_i$. Suppose for contradiction
that $D$ is  not $\kappa$-scattered.
As we may shrink $D$, without loss of generality we will assume that $D$ is strongly
$\kappa$-dense.
By induction on $\zeta<\kappa$, we define $i(\zeta), S_\zeta, T_\zeta$ as in
the following, if possible, and we stop at the first ordinal $\zeta^\ast$ for
which such a choice is not possible.

\underline{$i(0)$}. Let $i(0)$ be the first $i$ such that $R_i$ is of
size $\geq \kappa$,
which exists as $\kappa$ is regular.
As $R_{i(0)}$ is
$\kappa$-scattered, yet it has size $\geq \kappa$, we can find 
$S_0,T_0$ such that they are both subsets of $R_{i(0)}$ of size
$<\kappa$, and $S_0<_DT_0$, but for no $x\in R_{i(0)}$ do we have
$S_0<_Dx<_DT_0$.

\underline{$i(\zeta+1)$}. Given $i(\zeta), S_\zeta, T_\zeta$ such that
$S_\zeta$ and $T_\zeta$ are subsets of $R_{i(\zeta)}$ of size $<\kappa$ with
$S_\zeta<_D T_\zeta$, yet for no $x\in R_{i(\zeta)}$ do we have
$S_\zeta<_D x<_D T_\zeta$. As $D$ is strongly $\kappa$-dense
we have that $B_\zeta\deq\{x\in D:\,S_\zeta<_D x<_D T_\zeta\}$
has size at least $\kappa$ by Claim \ref{star}. By the regularity of $\kappa$, there is 
$i<i^\ast$ such that $B_\zeta\cap R_i$ has size at least $\kappa$. As 
$B_\zeta\cap R_{i(\zeta)} =\emptyset$, we have that first such $i$ is greater 
than $i(\zeta)$. We let this $i$ be $i(\zeta+1)$.

Now $B_\zeta\cap R_{i(\zeta+1)}$ is $\kappa$-scattered,
hence there are $S_{\zeta+1}$ and  $T_{\zeta+1}$ exemplifying this. In other
words, they are both subsets of $B_\zeta\cap R_{i(\zeta+1)}$ of
size $< \kappa$, and $S_{\zeta+1}<_D T_{\zeta+1}$, while for no $x\in B_\zeta
\cap R_{i(\zeta+1)}$  do we have $S_{\zeta+1} <_D x <_D T_{\zeta+1}$.

\underline{$i(\zeta)$ for $\zeta$ limit.} Let $i^+ \deq \mbox{ sup}_{\xi < \zeta} i(\xi)$. 
Hence, either $i^+ = i^*$, in which case we stop the induction, 
or $i^+ < i^*$, in which case we let $i(\zeta) \deq i^+$ and $S_{\zeta} \deq
\bigcup_{\xi < \zeta} S_{\xi}$, and similarly for $T_{\zeta}$. 

Notice that our induction must stop at some limit stage $\zeta^*  <
\kappa$  as $i^* < \kappa$. Now let $S \deq \bigcup_{\zeta <
\zeta^*} S_{\zeta}$, and similarly
for $T$. By the construction, it follows that $S <_D T$. Hence, there is $x \in
D$ with $S <_D x <_D T$. But then $x \in R_{i(\zeta + 1)}$ for some $\zeta$
(noticing that $D = \bigcup_{\zeta < \zeta^*} R_{i(\zeta + 1)}$  as the $R_i$'s
are increasing).  Therefore, $S_{\zeta + 1} <_D x <_D T_{\zeta + 1}$ and $x \in
R_{i(\zeta + 1)}$, a contradiction. 
\end{Proof of the Claim}

The analogue of the above claim is not
true for strongly $\kappa$-scattered posets, even when only the union of 
$\aleph_0$ strongly $\kappa$-scattered posets is considered; this follows from
the example in \S\ref{examples}, see Claim \ref{counterexample}.
However a weaker claim is true.

\begin{Claim}\label{katie} The union of two strongly $\kappa$-scattered
is strongly
$\kappa$-scattered. Consequently, the union of any finite number of such posets
is strongly $\kappa$-scattered.
\end{Claim} 

\begin{Proof of the Claim}
It clearly suffices to prove the first statement. So
assume that $D_1$ and $D_2$ are strongly $\kappa$-scattered
posets and assume that $D = D_1 \cup D_2$ embeds a weakly
$\kappa$-dense poset. By thinning out the sets, we can assume that
$D$ itself is weakly $\kappa$-dense. Then for all $s, t \in D$ we have
that the set $(s,t)_D = \{x \in D : s <_D x <_D t\}$ has size
$\geq \kappa$. 

On the other hand, there exist $s, t \in D_1$ such that $|(s,
t)_{D_1}| < \kappa$ as $D_1$ is strongly $\kappa$-scattered. Consider 
$T\deq\ (s,t)_{D}\rest D_2$.
As $T\subseteq D_2$, $T$ is not weakly $\kappa$-dense, and hence there
must be $u<_{D_2} v\in T $ such that $|(u,v)_{D_2}|<\kappa$.
But then $(u,v)_D=(u,v)_{D_2}\cup (u,v)_D\rest D_1$ must have size $<\kappa$,
as $(u,v)_{D}\rest D_1\subseteq (s,t)_{D_1}$. This is a contradiction.
\end{Proof of the Claim} 

Now we can prove the following lemma which holds both for weakly and
strongly $\kappa$-scattered posets. The version needed for the 
strongly $\kappa$-scattered case is to be read from within the square brackets.

\begin{Lemma}\label{equi} Assume that
$\kappa$ is a regular cardinal.
For any poset $P$, we have $(1)\implies (2)\implies (3)$:
\begin{description}
\item{(1)} $P$ is [strongly] $\kappa$-scattered and satisfies FAC,
\item{(2)} Every augmentation of $P$ is [strongly] $\kappa$-scattered,
\item{(3)} $P$ is [strongly] $\kappa$-scattered and satisfies
  $\lambda$-AC where $\lambda = \kappa^{< \kappa}$.
\end{description}
\end{Lemma}

\begin{Proof of the Lemma} Let $\kappa$ be as in the assumptions of the lemma.

$(1)\implies (2)$ Assume the contrary; $P$ is 
[strongly] $\kappa$-scattered and satisfies FAC, but $P'$ is an augmentation of $P$ 
which is not [strongly] $\kappa$-scattered. This implies 
that the size of $P$ and hence $P'$ is at least $\kappa$.

For the next subclaim, let  $(\perp q)^S$ be the set of all elements of $S$
that are incomparable with $q$ in $S$.

\begin{Subclaim}\label{sub2} Suppose that
$S \subseteq P$ and $A,B\subseteq S$ are such that $A<_S B$, 
$|A|,|B|<
\kappa$ [$|A|,|B|=1$] and 
$(A,B)_S=\emptyset$ [$|(A,B)_S|<\kappa$]
while with $S'=P'\rest S$ we have that 
$(A,B)_{S'}$ is not [strongly] $\kappa$-scattered [$|(A,B)_{S}'|=\kappa$].
Then there is $q \in A\cup B$ such that 
$(A,B)_{S'} \rest (\perp q)^S$ is not [strongly] $\kappa$-scattered.
\end{Subclaim}

\begin{Proof of the Subclaim}  Let $A, B$ and $S$ be as in the assumptions. 
 Let
\[\begin{array}{c} D_A = \{c \in S :  A <_{S'} c <_{S'} B \hbox{ and } a \perp_S c
\hbox{ for some }a\in A\}\\
D_B = \{c \in S : A <_{S'} c <_{S'} B \hbox{ and } b \perp_S c
\hbox{ for some }b\in B\}, \end{array}\]
so that $D_A \cup D_B = (A, B)_{S'}$.
Since the union of two [strongly] $\kappa$-scattered posets is itself
[strongly] $\kappa$-scattered by Claim \ref{mirna} [Claim \ref{katie}],
then either $S' \rest D_A$ or
$S' \rest D_B$ is not [strongly] $\kappa$-scattered.
[Let $q$ be the unique element of $A$ or of $B$, depending on which of the two
sets is not strongly $\kappa$-scattered, so finishing the proof in this case].

Now notice that 
\[
D_A= \bigcup_{a\in A}D_a\mbox{, where }D_a\deq\{c\in S:\,
A<_{S'}c<_{S'}B\mbox{ and }a\perp_S c\},
\]
and similarly for $D_B$.
Again by Claim \ref{mirna} there is
either $a\in A$ or $b\in B$ such that $D_a$ or $D_b$ is not
$\kappa$-scattered.  Let $q$ be either $a$ or $b$, whichever gives us the
non-$\kappa$-scattered poset.  Hence, $D_q=(A,B)_{S'} \rest (\perp q)^S$ is
not $\kappa$-scattered.   
\end{Proof of the Subclaim}

Suppose that $S \subseteq P$ is
such that $S'=P' \rest S$ is strongly [weakly]
$\kappa$-dense.  Such an $S$ exists
because $P'$ embeds a strongly [weakly] $\kappa$-dense
subset. Since $P\rest S $ is [strongly] $\kappa$-scattered
and has size $\kappa$ there must be
$A, B \subseteq S$ with $A <_P B$
and $|A|,|B|<\kappa$ with the property that the set $(A,B)_S$
is empty, as otherwise the condition $(\ast)^\kappa$ would be satisfied by
$S$ [there must be $A,B$ such that $|A|,|B|=1$ and 
$|(A,B)_S|<\kappa$].
By Claim \ref{star} we must have $|\{c \in S' : A
<_{S'} c <_{S'} B \mbox{ and }c\perp_S q\}| \geq \kappa$
[the analogue for strongly $\kappa$-scattered holds by the choice of $A,B$].
By induction on $n<\omega$ we shall choose
$A_n, B_n, q_n, S_n$ so that 
\begin{enumerate} 
\item $A_0=A$, $B_0=B$, $q_0=q$, $S_0=S$,
\item $|A_n|, |B_n|<\kappa$ [$A_n=\{a_n\}, B_n=\{b_n\}$], $A_n, B_n\subseteq S_n$,
\item $S_n$ is [strongly] $\kappa$-scattered while $S'_n=P'\rest S_n$ is not,
\item $A_n<_{S'_n} B_n$ and $(A_n, B_n)_{S_n}=\emptyset$ [$|(a_n,b_n)_{S_n}|<\kappa$ and
$|(a_n,b_n)_{S'_n}| \geq \kappa$],
\item $q_n\in A_n\cup B_n$ and $
(A_n, B_n)_{S'_n}\rest\bigcap_{k\le n}(\perp q_k)^{S_n}$
is not $\kappa$-scattered [$(a_n, b_n)_{S'} \rest (\perp q_n)_S$ is not strongly
$\kappa$-scattered],
\item $ A_{n+1}\cup B_{n+1}\subseteq (\perp q_n)^{S_n}$,
\item $A_n <_{S'_n} A_{n+1} <_{S'_n} B_{n+1} <_{S'_n} B_n$.
\end{enumerate} 
To start the induction we use the choices already made. Suppose that we are at the
stage $n+1$ of the induction. Since by the induction hypothesis $(A_n, B_n)_{S'_n}
\rest\bigcap_{\beta\le n}(\bot q_k)^{S_n}$ is not [strongly] $\kappa$-scattered,
it includes an order $T'_{n+1}$ which is strongly [weakly] $\kappa$-dense. Let
$T_{n+1}=P\rest T'_{n+1}$. Since $T_{n+1}$ is [strongly] $\kappa$-scattered
it does not satisfy $(\ast)^\kappa$ [it is not weakly
$\kappa$-dense] and hence there are $A_{n+1}, B_{n+1}
\subseteq T_{n+1}$ such that 
\[
|A_{n+1}|, |B_{n+1}|<\kappa \,\,\,[A_{n+1}=\{a_{n+1}\},
B_{n+1}=\{b_{n+1}\}]\mbox{ and }
A_{n+1}<_{T_{n+1}}B_{n+1}\]
but
\[
(A_{n+1}, B_{n+1})_{T_{n+1}}=\emptyset \,\,\,[|(a_{n+1}, b_{n+1})_{T_{n+1}}| <\kappa].
\]
Let $S_{n+1}=S_n
\cup T_{n+1}$ and hence $S'_{n+1}\supseteq S'_{n}\cup T'_{n+1}$
is $P'\rest S_{n+1}$.

We have $A_{n+1}, B_{n+1}\subseteq S_{n+1}$ and  $A_n<_{S'_n}
A_{n+1}<_{S'_n} B_{n+1} <_{S'_n} B_n$. Also
$S_{n+1}$ as the union of two [strongly] $\kappa$-scattered orders is
[strongly] $\kappa$-scattered by Claim \ref{mirna} [Claim \ref{katie}]
while $S'_{n+1}$ is not as it includes $T'_{n+1}$ which is strongly
[weakly] $\kappa$-dense. Note also that 
$A_{n+1}\cup B_{n+1}\subseteq \bigcap_{k\le n}(\bot q_k)^{S_n}$.

At any rate, Subclaim \ref{sub2} applies to $A_{n+1}, B_{n+1}$ and 
$S_{n+1}$ in place of $A,B$ and $S$. Hence we can find $q_{n+1}\in
A_{n+1}\cup B_{n+1}$ such that 
\[
(A_{n+1}, B_{n+1})_{S'_{n+1}}
\rest (\bot q_{n+1})^{S_{n+1}}=
(A_{n+1}, B_{n+1})_{S'_{n+1}}
\rest \bigcap_{k\le n+1}(\bot q_k)^{S_{n+1}}
\]
(as for $k\le n$ we have that $(A_{n+1}, B_{n+1})_{S'_{n+1}}
\subseteq (A_{n}, B_{n})_{S'_{n}} \subseteq 
(\bot q_k)^{S_{k}}$ $\subseteq (\bot q_k)^{S_{n+1}}$)
is not [strongly] $\kappa$-scattered, hence satisfying all the requirements of the induction at
this step.

Having finished the induction we obtain that if $k<n$ then $q_n\in A_n
\cup B_n$ $\subseteq (\bot q_k)^{S_{n}}\subseteq (\bot q_k)^{P}$,
hence $q_n\perp_P q_k$. Then the sequence $\langle q_n:\,
n<\omega\rangle$ forms an infinite antichain in $P$, contradicting the fact that
$P$ is FAC.

$(2)\implies (3)$
Suppose that every augmentation of $P$ is [strongly] $\kappa$-scattered but $P$ does
not satisfy the $\lambda$-AC for $\lambda = \kappa^{< \kappa}$.  ($P$ is automatically [strongly]
$\kappa$-scattered, since trivially $P$ is an augmentation of
itself.)  Take a subset $S \subseteq P$ such that $|S| = \lambda$ and
$S$ is a  $\lambda$-antichain.  We can now embed any strongly [weakly] $\kappa$-dense
set into $S$ so
forming an augmentation of $P$ which is not [strongly] $\kappa$-scattered. 
\end{Proof of the Lemma}

\begin{Remark} For $\kappa=\aleph_0$ the three conditions in Lemma \ref{equi}
are equivalent, as follows from the lemma. However for $\kappa>\aleph_0$ the
disjoint sum of an ordinal $\kappa$ with an antichain of size $\aleph_0$
shows that (3) does not imply (1) even for posets of size $\kappa$
when $\kappa = \kappa^{< \kappa}$. The above
proof does not seem to generalise to show that $(3)\implies (2)$ and we do not
know if this is the case.
\end{Remark}
 
\section{A Generalisation of the Classification}\label{mainsection}

Here we will generalise the classification of \cite{AB} to
$\kappa$-scattered FAC partial orders for regular $\kappa$.
From now on we will fix such a cardinal $\kappa$.
We remind the reader of the notion of the antichain rank of FAC posets, as introduced
in Definition \ref{ranks}.

\begin{Definition}  Fix some $\rho \geq 1$.  By induction on $\alpha$, an
ordinal, we define $\hal$ as follows:
\begin{description} \item{1. }  $\hnal_0 = \{1\}$.
\item{2. } $\hnal_1$ is the class of all posets $P$ satisfying the FAC
with $\hbox{ rk}_{\cal A}(P) \leq \rho$ such that either $P$ or its
inverse, or both, are $\kappa$-well founded.
\item{3. } If $\alpha$ is a limit ordinal, then 
$\hnal_\alpha = \bigcup_{\beta < \alpha} \hnal_\beta.$
\item{4. } If $\alpha = \beta + 1$ for some
$\beta>0$, then $\hnal_\alpha$ consists of all
posets $P$ that are lexicographical sums of the the form $P = \sum_{i \in I} P_i$
where $P_i \in \hnal_\beta$ and $I$ is in ${}^\kappa\HH^\rho_1$.
\end{description}

In general, let for all ordinals $\alpha$ and $\rho$,

\[
\hnal = \bigcup_{\alpha\mbox{ {\scriptsize{an ordinal}} }}\hal 
\hspace{1cm} \mbox{ and } \hspace{1cm}
\null^\kappa{\cal H} = \bigcup_{\rho\mbox{ {\scriptsize{an ordinal}} }} \hnal. 
\]
We let aug($\hnal$) be the set of all augmentations of posets in $\hnal$. 
\end{Definition}

\begin{Lemma}\label{oldhierarchy}
\begin{description} \item{(1)} The class $\hnal$ is the least class 
that contains the $\kappa$-well founded FAC posets with antichain ranks $\leq \rho$ 
and is closed under lexicographical sums and inverses. 
\item{(2)} Each $\hal$ and $\hnal$ is closed under restrictions and inverses.
\item{(3)} If $P \in \null^\kappa{\cal H}$ then $P$ is $\kappa$-scattered and satisfies the
FAC.
\item{(4)} aug$(\hnal)$ is closed under lexicographical sums, restrictions
and augmentations.  Every poset in aug$(\hnal)$ is $\kappa$-scattered.
\end{description}
\end{Lemma}

We remind the reader that $\kappa$-scattered is throughout used to refer to
weakly $\kappa$-scattered orders.

\begin{Proof of the Lemma}  {\bf (1)} It is clear that $\hnal$ contains all
$\kappa$-well founded posets of antichain ranks $\leq \rho$ and their
inverses, as $\hnal_1$ does. The proof that $\hnal$ is closed under  
lexicographical sums is the same as the one in \cite{AB} 
since we are holding $\kappa$ fixed. We will not include it here.  

It remains to show that $\hnal$ is closed under inverses. In fact we
shall prove 
by induction on $\alpha$ that each ${}^\kappa{\cal H}^\rho_\alpha$ is closed under inverses.
For $P \in \hnal$ we shall use the notation $\alpha(P) = \min\{\alpha :  P \in \hal\}$.
Let us commence the induction.

At $\alpha=0$ the situation is trivial and at
$\alpha = 1$, by definition $\hnal_1$ contains all inverses of its
members. 

At $\alpha = \beta + 1$, if $\alpha(P) < \alpha$ then this
case is covered by the induction hypothesis.  So, assume that 
$\alpha(P) = \alpha$. Then, $P = \sum_{i \in I}P_i$ where $P_i \in \hnal_{\beta}$ and 
$I \in \hnal_1$. The inverse of $P$ is $P^\ast = \sum_{i \in I^*} P_{i}^*$
where  $I^* \in \hnal_1$ because $\hnal_1$ is closed under inverses by definition, and
$P_{i}^*$ is the inverse of $P_i$.  We know that $P_i \in \hnal_{\beta}$,
thus $P_{i}^*$ is also in $\hnal_\beta$ by the induction hypothesis and hence $P^\ast$
is in $\hnal_\alpha$.

We know that $\alpha(P)$ is never a limit because $\alpha$ is a minimum.
Therefore, for $\alpha$ a limit ordinal and any $P\in
\hnal_\alpha$ , $\alpha(P)$ is strictly less than
$\alpha$.  Thus, this case is covered by the induction hypothesis. 

Hence $\hnal$ has the closure properties as required. Now we will show
that it is the least such class.  Suppose that
$\HH$ is another class with such properties. Again by induction on
$\alpha$, we will show that $\HH$ contains each $\hal$. Thus, we will show 
$\HH \supseteq \hnal$. The cases of $\alpha =0$ and $\alpha = 1$ are trivial by definition.
At $\alpha = \beta + 1$, all sets $P \in \hal$ are of the
form $P = \sum_{i \in I} P_i$ where each $P_i \in \hnal_\beta$.  Since
$\HH$ contains all $\hnal_\beta$ by the induction hypothesis and is closed
under lexicographical sums, all $P \in \hal$ must be in $\HH$.  Thus, 
$\hal \subseteq \HH$. The case where $\alpha$ a limit is similar since by definition, 
$\hal = \bigcup_{\beta < \alpha} \hnal_\beta$. 

{\bf (2)} We have already proved the closure under inverses in the proof of (1).
By induction on $\alpha$, we will show that each $\hal$ is closed under
restrictions. The case $\alpha=0$ is trivial.
For $\alpha = 1$, if $P \in \hnal_1$ then either $P$ or $P^\ast$ is $\kappa$-well founded.
Suppose that $P$ is $\kappa$-well founded.
Thus, if any restriction of $P$, call it $P^-$, had a $\kappa$-decreasing
sequence, it would actually be in $P$, which is a contradiction. The same
argument can be used for $P^*$, the inverse of any $\kappa$-well founded
poset in $\hnal_1$.

At $\alpha = \beta + 1$, suppose we are given $P = \sum_{i \in I} P_i$ where
each $P_i \in \hnal_{\beta}$ and $I\in \hnal_1$. By the induction hypothesis, all restrictions of
$P_i$ are in $\hnal_{\beta}$.  Any restriction, $P^-$, of $P$ can be
expressed as a lexicographical sum of restrictions of the $P_i$s along a restriction
of $I$. Thus
$P^-$ is also in $\hal$. The limit case is obvious.

{\bf (3)} Fix an ordinal $\rho$.
By induction on $\alpha$, we will prove that any $P \in \hnal_\alpha$ is
$\kappa$-scattered. The case $\alpha=0$
is trivial. For $\alpha = 1$, notice that since any strongly
$\kappa$-scattered order has a
$\kappa$-decreasing sequence by Lemma \ref{q-dec}, we have that  no $\kappa$-well
founded poset could embed such an order. Similarly, since by the same lemma
strongly $\kappa$-dense orders have $\kappa$-increasing sequences, a poset whose inverse is
$\kappa$-well founded also cannot embed such an order.
The limit case of the induction is taken care of by the induction hypothesis.

For $\alpha = \beta + 1$, if $P \in \hal$ we can
by the induction hypothesis let 
$P = \sum_{i \in I} P_i$ where each $P_i$ is $\kappa$-scattered and $I$ is
$\kappa$-scattered.  We will show that $P$ is
$\kappa$-scattered. For the sake of contradiction, let $Q$ be a
strongly $\kappa$-dense order and suppose $f: Q
\rightarrow P$ is an order preserving embedding. 

\underline{Case 1}. For every $i \in I$, there is at most one $q \in Q$ such
that $f(q) \in P_i$.  Define $g: Q \rightarrow I$ by letting $g(q)
= i$ iff $f(q) \in P_i$. This is well-defined by the assumptions of Case 1.
We also have $q <^* r$ implies $f(q) <_P f(r)$  which
implies $g(q) <_I g(r)$
by the definition of the lexicographic sum. Hence, $g$ is an order preserving embedding,
contradicting the fact that $I$ is $\kappa$-scattered.  

\underline{Case 2}. Not Case 1. There is an $i \in I$  and $q, r \in Q$ such
that $q \neq r$ and $f(q), f(r) \in P_i$.  Without loss of generality, take
$q <^\ast r$.  Because $f$ is an embedding, $f(q) <_{P_i} f(r)$.  By the
definition of the sum, we also have $f(x) \in P_i$ for all $x \in (q, r)$.
However, $(q, r)$ is strongly $\kappa$-dense by Claim \ref{star}, so $P_i$ is not
$\kappa$-scattered, which is a contradiction.

A similar proof shows that the second part of the claim in (3) is true.

{\bf (4)} The second sentence has already been covered in Lemma
\ref{equi}, because (3) shows that every element of
$\hnal$ satisfies the statement (1) of that lemma.
The first sentence of (4) is easily proven with each property
requiring the same type of argument. For example, to prove that
aug($\hnal$) is closed under lexicographical sums, consider the
following observations: suppose $P=\sum_{i\in J} P_i$ and each $P_i$ is an
augmentation of $Q_i$, where $Q_i$ is in $\hnal$ while $J$ is
an augmentation of $I\in \hnal$. Then $P$ is an
augmentation of $Q=\sum_{i\in I} Q_i$. As we have that $Q$ belongs to $\hnal$
(by part (1)), we conclude that $P$ is in aug($\hnal$).
\end{Proof of the Lemma}

The next theorem is virtually the same in claim and proof as Theorem
2.3 of
\cite{AB}. The only modification is the larger classification $\hal$ that
replaces the $\HH^{\rho}_\alpha$ in the paper. The proof for the
larger classification is the same because we hold $\kappa$ fixed, as
we have done with all other proofs of this nature. We will leave this
theorem as a fact, rather than reiterating the proof. 

Before we state the theorem, we need to draw attention to an unusual ordinal 
operation known as Hessenberg based exponentiation. This smoothly extends the
Hessenberg product operation which in turn
extends the natural sum operation. Since we do not 
need to know the exact value of the exponent for this paper, we refer
the reader to \cite{AB} for a more precise definition. 
We denote the Hessenberg 
based exponentiation of $\alpha$ and $\beta$ by $\alpha^{H\beta}$.   

\begin{Theorem} If $P \in \hal$ then rk$_{{\cal A}}(P) \leq \rho^{H\alpha}$.
\end{Theorem}

Hausdorff's theorem \cite{Ha} (or see \cite{haus})
and the Abraham-Bonnet generalisation in \cite{AB} are both {\em characterisations}
of the class of linear and FAC posets, respectively, which do not embed the rationals.
The latter class is exactly ${}^{\aleph_0}{\cal H}$.
To prove something like that we would need to know that
if $P \not\in \mbox{aug}({}^\kappa{\cal H})$ is an FAC poset, then $\qkappa$
embeds into $P$.  Unfortunately we have not been able to prove such a claim
for uncountable $\kappa$, and the question if it
is true even if we assume that
$\kappa$ has some large cardinal properties
remains open. We shall instead prove a weaker claim, for which we shall fatten up our hierarchy 
a little.

\begin{Definition}\label{Hstar} Let ${}^\kappa{\cal H}^\ast$ denote the closure of
\mbox{aug}(${}^\kappa{\cal H})$ under FAC weakenings, that is, the class obtained by taking all
FAC orders $P$ for which there is an order $P'$ in ${}^\kappa{\cal H}$ such that $P$
is a weakening of $P'$.
\end{Definition}

We shall show that ${}^\kappa{\cal H}^\ast$ lies between the classes of 
strongly and
weakly $\kappa$-scattered FAC partial orders. Let us first show the easy direction.

\begin{Claim}\label{up} Every poset in ${}^\kappa{\cal H}^\ast$ is (weakly)
$\kappa$-scattered and FAC.
\end{Claim}

\begin{Proof of the Claim}
Let $P$ be in ${}^\kappa{\cal H}^\ast$ and let $P'$ in $\mbox{aug}({}^\kappa{\cal H})$
be such that $P$
is an FAC weakening of $P'$. Clearly $P$ is FAC. If $P$ were not to be 
weakly $\kappa$-scattered then some strongly $\kappa$-dense order
would embed into $P$ and hence into $P'$, in 
contradiction with Lemma \ref{oldhierarchy}(3) and Lemma \ref{oldhierarchy}(4).
\end{Proof of the Claim}

The heart of our main theorem lies in the following:

\begin{Claim}\label{includesall}
Every strongly $\kappa$-scattered FAC partial order belongs to
${}^\kappa{\cal H}^\ast$.
\end{Claim}

\begin{Proof of the Claim} 
Suppose for contradiction that $P$ is a strongly $\kappa$-scattered FAC partial order
which does not belong to
${}^\kappa{\cal H}^\ast$. Recalling that every poset has a linear augmentation, 
let $Q$ be any linear augmentation of $P$. By Lemma \ref{equi}
$Q$ is strongly $\kappa$-scattered, and by the definition of ${}^\kappa{\cal H}^\ast$
we have that $Q\notin {}^\kappa{\cal H}$ (and even
$Q\notin \mbox{aug}({}^\kappa{\cal H})$).

For $a,b\in Q$ we define an equivalence relation $a\equiv b$ iff
the interval in $Q$ between $a$ and $b$ is in ${}^\kappa{\cal H}$. It is easily seen that this
indeed is an equivalence relation. For $a\in Q$ let $C_a=\{b:\,a\equiv b\}$ be
the equivalence class of $a$.

\begin{Subclaim}\label{first} Each $C_a$ with the order induced from $Q$ is in ${}^\kappa{\cal H}$.
\end{Subclaim}

\begin{Proof of the Subclaim} Given $C_a$. By induction on $\gamma$ an ordinal pick
if possible $a_\gamma$ and $b_\gamma$ in $C_a$
so that $a_0=b_0=a$, $a_\gamma$ is $Q$-increasing
with $\gamma$ and $b_\gamma$ is $Q$-decreasing with $\gamma$. Since $C_a$ is a set,
there must be ordinals $\alpha$, the first $\gamma$ for which we cannot choose
$a_\gamma$ and $\beta$, the first $\gamma$ for which we cannot choose
$b_\gamma$. Then $C_a$ is the lexicographic sum
\[
\Sigma_{i<\beta}[b_{i+1}, b_i)\oplus\{a_0\}\oplus\Sigma_{j<\alpha}(a_j,a_{j+1}].
\]
Note that each of the intervals mentioned above is ${}^\kappa{\cal H}$, by the definition of
$C_a$ and the fact that $\equiv$ is an equivalence relation.
Since ${}^\kappa{\cal H}$ is closed under lexicographic sums of the above kind, we obtain
that $C_a\in {}^\kappa{\cal H}$.
\end{Proof of the Subclaim}

\begin{Subclaim}\label{second} If $a, b \in Q$ are not
$\equiv$-equivalent, and $a<_Q b$, then $C_a <_Q C_b$.
\end{Subclaim}

\begin{Proof of the Subclaim} Let $c\in C_a$ and $d\in C_b$. 
Clearly $c\neq d$. Suppose for contradiction that $d<_Q c$ and distinguish
two cases.

\underline{Case 1}. $c\le _Q a$.

Then $d<_Q a$, so $(a,b)\subseteq (d,b)$, which is a contradiction because the latter
is a member of ${}^\kappa{\cal H}$ while the former is not.

\underline{Case 2}. $c\ge _Q a$.

Then either $d<_Q a$, in which case we obtain a contradiction like in Case 1, or $a\le_Q d$.
In the latter case we have that $(d,c)\subseteq (a,c)$, contradicting the fact that
$(a,c)\in {}^\kappa{\cal H}$ and $(d,c)\notin {}^\kappa{\cal H}$.
\end{Proof of the Subclaim}

Let $Q^\ast$ be a set of the representatives of the $\equiv$-equivalence classes ordered
by the factor order (by Subclaim \ref{second} this order agrees with the order in $Q$).
Then $Q$ is the lexicographic sum $\Sigma_{a\in Q^\ast} C_a$ and since $Q\notin {}^\kappa{\cal H}$
we obtain by Subclaim \ref{first} and the closure of ${}^\kappa{\cal H}$
 under lexicographic sums
that $Q^\ast\notin {}^\kappa{\cal H}$. In particular $Q^\ast$ has size at least $\kappa$.

Now note that by the choice of $Q^\ast$ for every $a <_Q b$ in $Q^\ast$ the interval
$(a,b)_Q$ is not in ${}^\kappa{\cal H}$
(and that $Q^\ast$ is a maximal such set). We claim that
in fact $(a,b)_{Q^\ast}\notin {}^\kappa{\cal H}$ for such $a,b$. Once we prove this we shall be done,
because every poset of size $<\kappa$ is easily seen to be in ${}^\kappa{\cal H}$ and thus, 
$Q^*$ is a weakly $\kappa$-dense subset of $Q$.

So suppose that $a<_Q b$ are elements of $Q^\ast$ but $(a,b)_{Q^\ast}\in {}^\kappa{\cal H}$.
We then observe that $(a,b)_Q$ is the lexicographic sum $\Sigma_{c\in (a,b)_{Q^\ast}}
C_c$, which would then have to be in ${}^\kappa{\cal H}$, a contradiction.
\end{Proof of the Claim}

To finish our work we
shall give a simpler description of the class ${}^\kappa{\mathcal H}^\ast$.
We show that we do not need to start with $\kappa$-well-founded FAC posets
in the formation of ${}^\kappa{\mathcal H}^1$, we may start with
$\kappa$-well founded linear orders and then the FAC posets get picked up
when we form ${}^\kappa{\mathcal H}^\ast$. 

\begin{Claim}\label{simplerHast} Suppose that $\kappa$ is a regular cardinal.
Then ${}^\kappa{\mathcal H}^\ast$
is the closure of the class of all $\kappa$-well founded linear orders under
inversions, lexicographic sums, FAC weakenings and augmentations.
\end{Claim}

\begin{Proof of the Claim} Let $\mathcal H$ denote the closure of of the class 
of all $\kappa$-well founded linear orders under
inversions, lexicographic sums, FAC weakenings and augmentations. Since
${}^\kappa{\mathcal H}^\ast$ is the closure of the class of 
$\kappa$-well founded FAC posets under these operations we have that
${}^\kappa{\mathcal H}^\ast\supseteq H$. On the other hand, if $P\in 
{}^\kappa{\mathcal H}^\ast$ then let $Q\in {\rm aug}({}^\kappa{\mathcal H})$ be
such that $P$ is an FAC weakening of $Q$ and let $R\in {}^\kappa{\mathcal H}$ be
such that $Q$ is an augmentation of $R$. If $R\in {\mathcal H}$ then 
$Q\in {\mathcal H}$ by the closure of ${\mathcal H}$ under augmentations and hence
$P\in {\mathcal H}$ by the closure of ${\mathcal H}$ under FAC weakenings.
Hence it suffices to show that ${}^\kappa{\mathcal H}\subseteq {\mathcal H}$.
Let $\rho\ge 1$ be any ordinal, we shall show 
by induction on $\alpha$ that ${}^\kappa {\mathcal H}^\rho_\alpha
\subseteq {\mathcal H}$. We first need a subclaim.

\begin{Subclaim}\label{largecard} 
Every augmentation of a $\kappa$-well founded FAC poset is
$\kappa$-well founded.
\end{Subclaim}

\begin{Proof of the Subclaim} Let $P$ be a $\kappa$-well founded FAC poset
and $Q$ an augmentation of $P$. Suppose that $\langle a_\alpha;
\,\alpha<\kappa\rangle$ is a $\le_Q$-decreasing sequence. For
$\alpha<\beta<\kappa$ define $f(\alpha,\beta)=1$ if $a_\alpha$ and
$a_\beta$ are comparable in $P$ and let $f(\alpha,\beta)=0$ otherwise.
We now use the Dushnik-Miller theorem which says that either there is
an infinite 0-homogeneous set or a 1-homogeneous set of type $\kappa$.
Since $P$ is an FAC poset there cannot be an infinite 0-homogeneous
set, but a 1-homogeneous set of type $\kappa$ would contradict
the fact that $P$ is $\kappa$-well founded. This contradiction proves the subclaim.
\end{Proof of the Subclaim}

We now proceed with the promised inductive proof.
If $\alpha=0$ the conclusion is clear. If $P\in
{}^\kappa{\mathcal H}^\rho_1$ then $P$ is FAC and either $P$ or its
inverse (or both) are $\kappa$-well founded.

In the first case we can use the subclaim to find $Q$
which is a $\kappa$-well founded linear augmentation of $P$. Hence
$Q\in {\mathcal H}$ and as its FAC weakening, $P\in {\mathcal H}$.
The other case is similar.

The case of $\alpha$ a limit ordinal follows from the inductive hypothesis
and the case $\alpha=\beta+1$ for $\beta>0$ follows by the closure of
${\mathcal H}$ under lexicographic sums.
\end{Proof of the Claim}

Let us also observe the following:

\begin{Observation} Suppose that $P$ is a linear order, $Q$ is an FAC
  weakening of $P$,
and $R$ is an augmentation of $Q$. Then $R$ is an FAC weakening of $P$.
\end{Observation}

We conclude that the following theorem is true.

\begin{Main Theorem}\label{main}
Assume that $\kappa$ is a regular cardinal. Let
${}^\kappa\HH^\ast$ denote the closure of the class
of all $\kappa$-well founded linear orders
under inversions, lexicographic sums and FAC weakenings. Equivalently,

{\noindent (1)} 
${}^\kappa\HH^\ast$ 
contains all strongly $\kappa$-scattered FAC posets. 

{\noindent (2)} ${}^\kappa\HH^\ast$ is contained in the class of all $\kappa$-scattered
FAC posets.
\end{Main Theorem}

If $\kappa=\aleph_0$ we obtain an equality between the notions of $\kappa$-dense
and strongly $\kappa$-dense, so
applying Theorem \ref{main} to $\kappa=\aleph_0$ we obtain
that $\HH^\ast$ is the class of all scattered FAC posets. Since Abraham-Bonnet
theorem already gives that this class of posets is described by $\HH$ we have as
a corollary

\begin{Corollary} ${}^{\aleph_0}\HH^\ast$ is exactly the Abraham-Bonnet class
${}^{\aleph_0}\HH$.
\end{Corollary}

In general the two notions of density are not equivalent, as we illustrate in
\S\ref{examples}. Moreover, example \ref{counterexample} shows that for every
uncountable $\kappa$ with $\kappa=\kappa^{<\kappa}$ there are members of
${}^\kappa{\cal H}$ which are not strongly $\kappa$-scattered.
We also do not know for which uncountable $\kappa$ we obtain that ${}^{\kappa}\HH^\ast$
is the same as
${}^{\kappa}\HH$. Note that it is not to be expected that ${}^\kappa\HH^\rho$ is closed
under FAC weakenings as weakening a partial order generally adds larger antichains and hence
increases the antichain rank.

When reduced to the class of linear orders the class ${}^{\kappa}\HH^\ast$ can be 
replaced by a simpler class.

\begin{Theorem}\label{classL}
Assume that $\kappa$ is a regular cardinal. Let ${}^\kappa{\mathcal L}^\ast$
denote the closure of the class of all $\kappa$-well founded linear orders under inversions
and lexicographic sums. Then:

{\noindent (1)} ${}^\kappa{\mathcal L}^\ast$ contains all strongly $\kappa$-scattered
linear orders.

{\noindent (2)} ${}^\kappa{\mathcal L}^\ast$ is contained in the class of all 
$\kappa$-scattered linear orders.
\end{Theorem}

\begin{Proof} Linear orders are FAC posets with antichain rank $\le 1$. By Lemma 
\ref{oldhierarchy}(1)
the class  ${}^\kappa{\mathcal H}^1$ is the least class that contains the
$\kappa$-well founded linear orders and is closed under inversions and lexicographic
sums, hence ${}^\kappa{\mathcal L}^\ast={}^\kappa{\mathcal H}^1$. Since every
order in ${}^\kappa{\mathcal H}^1$ is linear we obtain ${}^\kappa{\mathcal H}^1=
{\rm aug}({}^\kappa{\mathcal H}^1)$, and hence Lemma 
\ref{oldhierarchy}(4) gives part (2) of the theorem.

To prove (1) we use the proof of Claim \ref{includesall}. We start with a strongly
$\kappa$-scattered linear order $Q$ that does not belong to ${}^\kappa{\mathcal L}^\ast
={}^\kappa{\mathcal H}^1$ and obtain a contradiction literally as in the proof of that
claim.
\end{Proof}

With $\kappa=\aleph_0$ Theorem \ref{classL} gives Hausdorff's theorem.

The above theorems and remarks raise the following questions

\begin{Question}\label{open} {\noindent (1)} For which uncountable $\kappa$ is
aug($\hnal$) exactly the class of all $\kappa$-scattered
FAC posets with antichain rank $\le\rho$? 

{\noindent (2)} For which $\kappa$ is it true that any FAC poset all
of whose
subposets (or even just chains) of size $\kappa$ belong to ${}^\kappa\HH$, is
itself an element of ${}^\kappa\HH$?

{\noindent (3)} For which $\kappa>\aleph_0$ is it true that every augmentation of a
[strongly] $\kappa$-scattered $\kappa$-AC poset is [strongly] $\kappa$-scattered?

{\noindent (4)} For which $\kappa>\aleph_0$ is ${}^\kappa\HH$ closed under FAC weakenings?
\end{Question}

We comment that one may generalise Theorem \ref{main} to the case of $\lambda
<\kappa$ where both $\lambda$ and $\kappa$ are equal to their weak powers, and
consider the situation of posets of size $\kappa$ that satisfy (strong) $\lambda$-density,
obtaining the expected results.

\section{Examples}\label{examples}

For the sake of completeness we include some examples that illustrate the difference
between weak $\kappa$-density and strong $\kappa$-density. We shall assume that 
$\kappa$ is an uncountable regular cardinal.

An easy example of a linear order that is $\kappa$-dense but not strongly $\kappa$-dense
is the lexicographic sum along $\omega+\omega^\ast$ of any strongly
$\kappa$-dense order. This order is clearly strongly
$\kappa$-dense. We give an example of a $\kappa$-dense linear order 
which is weakly $\kappa$-scattered and moreover does not have a
$\kappa$-decreasing sequence. 

Let $L_0$ be the lexicographic sum along $\omega^\ast$ of copies of $\kappa$.
By induction on $n<\omega$ define $L_n$ by letting $L_{n+1}$ be
the lexicographic sum along $L_n$ of copies of $L_0$.
We denote the order of $L_n$ by $\le_n$.
Let $L=\bigcup_{n<\omega} L_n$ be
ordered by letting $p\le q$ iff $p\le_n q$ for the first $n$ that contains both $p$
and $q$.

\begin{Claim}\label{counterexample}
No $L_n$ for $n<\omega$ is $\kappa$-dense. $L$ is $\kappa$-dense.
\end{Claim}

\begin{Proof} The first statements can easily be proven by induction. For the
second one,
let $p < q$ and let $n$ be the first such that $p,q\in L_n$.
By the definition of $L_{n+1}$ there is a copy of $L_0$ in $\{x\in
L_{n+1}:\,p<x<q\}$, so clearly the size of this set is $\kappa$.
\end{Proof}

\begin{Claim} $L$ does not have a decreasing sequence of size $\kappa$.
\end{Claim}

\begin{Proof} Suppose it had such a decreasing sequence, call it $S$. Then
$S=$ $\bigcup_{n<\omega} S\cap (L_{n+1}\setminus L_n)$.
By the regularity of $\kappa>\aleph_0$
there has to be $n$ for which the size of $S\cap (L_{n+1}\setminus L_n)$ is
$\kappa$, hence it suffices for us to show that no $L_n$ can have a decreasing
sequence of size $\kappa$. This can be done by induction on $n$.
\end{Proof}

Hence by Lemma \ref{q-dec} we have 

\begin{Corollary} $L$ does not embed any strongly $\kappa$-dense order
  and $L$ is in $\hnal$ for any $\rho\ge 1$.
\end{Corollary}

This shows that the
boundary of $\hnal$ is somewhere in between weakly $\kappa$-scattered
and strongly $\kappa$-scattered. We
conjecture that  $\hnal$ contains all strongly $\kappa$-scattered FAC
posets.

\bibliographystyle{plain}
\bibliography{shortref}

\begin{thebibliography}{1}

\bibitem{AB}
U.~Abraham and R.~Bonnet.
\newblock Hausdorff's theorem for posets that satisfy the finite antichain
  property.
\newblock {\em Fundamenta Mathematica}, 159(1):51 -- 69, 1999.

\bibitem{haus}
G.~Asser, J.~Flachsmeyer, and W.~Rinow.
\newblock {\em Theory of Sets and Topology; In honour of Felix Hausdorff}.
\newblock Deutscher Verlag der Wissenschaften, 1972.

\bibitem{BP}
R.~Bonnet and M.~Pouzet.
\newblock Linear extensions of ordered sets.
\newblock In {\em Ordered Sets}, pages 125 -- 170. D. Reidel Publishing
  Company, 1982.

\bibitem{fraisse}
R.~Fra{\"{\i}}ss{\'{e}}.
\newblock {\em Theory of Relations}, volume 145 of {\em Studies in Logic and
  Foundations of Mathematics}.
\newblock Elsevier Science, B.V., revised edition, 2000.

\bibitem{Ha}
F.~Hausdorff.
\newblock {Grundz\"uge} einer {T}heorie der geordnete {M}engenlehre.
\newblock {\em Mathematische Annalen}, 65:435 -- 505, 1908.
\newblock (in German).

\bibitem{linorders}
J.~Rosenstein.
\newblock {\em Linear Orderings}.
\newblock Pure and Applied Mathematics. Academic Press, 1982.

\bibitem{Shbookc}
S.~Shelah.
\newblock {\em Classification Theory}, volume~92 of {\em Studies in Logic and
  Foundations of Mathematics}.
\newblock North-Holland, revised edition, 1990.

\end{thebibliography}

\bigskip

\parbox{7cm}{
Corresponding author: \\
Mirna D\v zamonja\\
h020@uea.ac.uk\\
School of Mathematics\\
University of East Anglia \\
Norwich, NR4 7TJ, UK\\
}

\parbox[b]{7cm}{
Katherine Thompson\\
aleph\_nought@yahoo.com\\
Department of Mathematical Sciences\\
Carnegie Mellon University \\
Pittsburgh, PA 15213 USA\\}

\end{document}